\documentclass[12pt]{article}
\usepackage[OT2,T1]{fontenc}
\DeclareSymbolFont{cyrletters}{OT2}{wncyr}{m}{n}
\DeclareMathSymbol{\Sha}{\mathalpha}{cyrletters}{"58}
\usepackage{amsmath,amssymb,url}
\usepackage[none]{hyphenat}
\begin{document}
{\Large\begin{center}{\bf Large Tate--Shafarevich orders from good $abc$ triples}\\
David Broadhurst, Open University, UK\\
{\tt David.Broadhurst@open.ac.uk}\\
14 November 2021
\end{center}}{\large
Record values are determined for 
the order $|\Sha|$ of the Tate--Shafarevich group of an
elliptic curve $E$, computed analytically by the
Birch--Swinnerton-Dyer conjecture,
and for the Goldfeld--Szpiro ratio
$G=|\Sha|/\sqrt{N}$, where $N$ is the conductor of $E$. 
The curves have rank zero and are isogenous to 
quadratic twists of Frey curves constructed from
coprime positive integers $(a,b,c)$ with $a+b=c$
and $c>r^{1.4}$, where the radical $r$ is the 
product of the primes dividing $abc$. 
Curves with  $|\Sha|>250000^2$ and $G>12$ are found in 20 isogeny classes.
Three curves have $G>150$. The largest value of
$|\Sha|$ is $1937832^2>3.755\times10^{12}$. This is
more than 3.5 times the previous record, which had been
computed at a cost about 600 times greater than that for the new record. The primes
25913, 27457, 36929 and 49253 are identified as divisors of $|\Sha|$ values.

\newpage 

\section{Definitions and results}

Let $a$ and $b$ be coprime integers with $b>a>0$. Let $c=a+b$ and
$\lambda=\log(c)/\log(r)$, where the radical $r$ is the product of the primes $p|abc$.
There are 241 known $abc$ triples~\cite{S} with $\lambda>1.4$. Given such a triple
and a square-free integer $q$, one may construct an isogeny class
of 4 elliptic curves over the rationals~\cite{W}, as follows:
\begin{eqnarray}
E_1&:&y^2=x\left(x^2-2q(b+c)x+q^2a^2\right)\label{E1}\\
E_2&:&y^2=x\left(x^2+2q(a+c)x+q^2b^2\right)\label{E2}\\
E_3&:&y^2=x\left(x^2-2q(a-b)x+q^2c^2\right)\label{E3}\\
E_4&:&y^2=x(x+qb)(x+qc)\label{E4}
\end{eqnarray}
where $E_4$ is a quadratic twist of a Frey curve. If these curves have rank zero,
then the conjecture of Birch and Swinnerton-Dyer asserts that the
central value of the $L$-function associated with this isogeny class is given by
\begin{equation}
L=L(1)=\frac{\pi C_k|\Sha(E_k)|}{\sqrt{|q|c}\,{\rm AGM}(1,\alpha)},\quad
\alpha=\left\{\begin{array}{l}\sqrt{b/c}~{\rm for}~q>0\\\sqrt{a/c}~{\rm for}~q<0
\end{array}\right.\label{agm}
\end{equation}
where $C_k$ is a positive integer, determined by Tamagawa factors at primes
dividing the conductor $N$, and
$|\Sha(E_k)|$ is the {\em order} of the Tate-Shafarevich group $\Sha$ of $E_k(c,a,q)$
and is assumed to be the square of an integer. The integers $C_k$ 
differ only by powers of 4. Tables 1 to 4 give the smallest value
of $k$ for which $C_k$ is minimal and
the corresponding Goldfeld-Szpiro ratio~\cite{DS,DW,GS} $G=|\Sha|/\sqrt{N}$ is hence maximal.

\newpage

{\bf Table 1}: Curves with $|\Sha|>250000^2$ and $G>12$.
\begin{gather*}\begin{array}{|c|c|c|c|c|c|c|}\hline|\Sha|&c&a&q&k&L&G\\\hline
1937832^2&3^{19}11^{4}463^{5}&5^{4}19^{13}103&285&3&21.1540&153.084\\
804572^2&3^{22}7^{14}43\cdot83&2^{5}67^{8}107\cdot22381&7&2&5.64497&163.119\\
793656^2&2^{9}5^{16}11^{9}79&29^{4}2213^{2}&23&1&17.3059&162.256\\
589080^2&3^{19}11^{4}463^{5}&5^{4}19^{13}103&-95&1&6.98206&24.5023\\
574656^2&3^{13}5^{8}11^{3}53\cdot73^{2}89^{2}103&7^{5}61&39&1&4.61571&24.4672\\
514672^2&11^{8}109^{2}3677^{3}&2\cdot5^{10}13^{4}&429&3&7.22644&24.4535\\
487408^2&2\cdot3^{15}7^{2}31^{10}&5\cdot67^{3}127^{2}19219&62&3&2.82526&20.0273\\
480512^2&2^{37}89^{3}167^{2}1823&3^{22}9787^{2}&29&1&4.85512&28.8220\\
479144^2&5^{9}139^{6}&2^{2}11&104945&1&37.9640&96.2939\\
439312^2&5^{9}139^{6}&2^{2}11&114998&1&30.4875&54.6806\\
421216^2&5^{11}7^{10}79\cdot389^{2}&11\cdot103^{8}&778&1&0.78741&12.7045\\
394024^2&2\cdot3^{15}7^{2}31^{10}&5\cdot67^{3}127^{2}19219&26&3&2.85120&20.2112\\
393216^2&3^{38}13^{4}5233&71^{8}233^{3}&5&1&1.46770&34.8223\\
338122^2&3^{38}397&13^{5}19^{3}&5161&1&20.7254&16.9191\\
324440^2&5^{15}179^{4}2141&2^{12}13^{3}223^{3}&30&3&11.1895&39.8368\\
321584^2&5^{9}139^{6}&2^{2}11&5282&3&38.1136&31.3650\\
298500^2&2\cdot11^{6}193^{4}20551&3\cdot5^{6}7^{8}53&5010&3&25.1811&34.9472\\
295432^2&3^{13}5^{8}11^{3}53\cdot73^{2}89^{2}103&7^{5}61&1&1&11.4278&14.2781\\
288182^2&2^{4}3^{19}17^{8}29&5^{2}23^{10}106531&108273&2&13.3365&15.3808\\
277548^2&3^{30}13^{4}277&5^{11}31\cdot191&4966&3&5.44574&14.7707\\
\hline\end{array}\end{gather*}

\newpage

{\bf Table 2}: More curves with $|\Sha|>10^{10}$ and $G>7$.
\begin{gather*}\begin{array}{|c|c|c|c|c|c|c|}\hline|\Sha|&c&a&q&k&L&G\\\hline
245520^2&2\cdot11^{6}193^{4}20551&3\cdot5^{6}7^{8}53&10615&2&5.85180&16.2426\\
243952^2&11^{8}109^{2}3677^{3}&2\cdot5^{10}13^{4}&186&1&4.93145&8.34373\\
240848^2&5^{9}139^{6}&2^{2}11&807&3&27.3470&23.9241\\
217564^2&2^{15}17^{2}331\cdot1061^{4}&3^{3}241^{3}&10065&1&12.0412&20.1124\\
211784^2&5^{9}139^{6}&2^{2}11&272857&1&9.19960&11.6672\\
205720^2&5^{15}17\cdot53093^{2}&2^{11}3^{4}101^{4}29221&5151&2&6.20109&18.2946\\
202760^2&2^{26}5^{12}1873&13^{10}37^{2}&49913&1&6.67897&22.8167\\
197256^2&2\cdot3^{3}5^{23}953&7^{2}41^{2}311^{3}&63755&1&5.00374&7.17519\\
185028^2&5^{9}11^{8}2489197589&3^{4}7^{2}41&33&1&6.16228&25.7619\\
183200^2&5^{9}139^{6}&2^{2}11&92883&3&14.7483&7.99828\\
160368^2&3^{30}13^{4}277&5^{11}31\cdot191&4011&1&4.04598&7.75986\\
127542^2&11^{9}29^{4}101^{3}&2^{2}3^{4}163^{3}1006151&163&1&6.18013&12.2005\\
116224^2&11^{12}389^{2}6841&2^{6}5^{2}7^{13}13^{2}463&-1945&2&5.75633&7.67883\\
114704^2&2\cdot3^{3}5^{23}953&7^{2}41^{2}311^{3}&-1555&2&9.52652&15.5354\\
112704^2&2\cdot3^{3}5^{23}953&7^{2}41^{2}311^{3}&4043&1&0.81083&9.30159\\
110268^2&2\cdot3^{3}5^{23}953&7^{2}41^{2}311^{3}&3157&1&5.27005&10.0761\\
109976^2&2\cdot11^{6}193^{4}20551&3\cdot5^{6}7^{8}53&1837&2&8.46717&7.83401\\
106064^2&2^{30}5^{2}127\cdot353^{2}&1&10590&1&8.43066&14.0020\\
103312^2&5^{15}17\cdot53093^{2}&2^{11}3^{4}101^{4}29221&-39&2&7.68388&18.7474\\
103296^2&2^{26}5^{12}1873&13^{10}37^{2}&9365&1&10.0047&13.6712\\
\hline\end{array}\end{gather*}

\newpage

{\bf Table 3}: More curves with $|\Sha|>10^{8}$ and $G>12$.
\begin{gather*}\begin{array}{|c|c|c|c|c|c|c|}\hline|\Sha|&c&a&q&k&L&G\\\hline
87920^2&3^{9}5^{6}13^{5}23\cdot191&17^{4}&4485&1&16.3832&13.7675\\
83168^2&5^{9}139^{6}&2^{2}11&1390&1&4.96930&17.8253\\
70400^2&5^{9}139^{6}&2^{2}11&302&1&7.63893&27.4015\\
63036^2&3^{38}397&13^{5}19^{3}&1&1&12.9372&42.2448\\
56932^2&2^{28}3^{12}11^{3}67&7^{11}19&70&1&10.9192&20.9225\\
40800^2&19^{3}139^{6}&3\cdot13^{2}1049&278&3&45.1256&30.4809\\
40800^2&2^{11}73^{7}83^{2}197&3^{5}5^{15}13^{5}&1&1&3.34266&17.4567\\
34168^2&19^{11}59\cdot7207&2^{46}23&6&3&3.40399&13.1792\\
33600^2&3^{16}103^{3}127&73&2823&1&6.90932&38.6702\\
32272^2&19^{3}139^{6}&3\cdot13^{2}1049&1807&3&5.53691&29.9201\\
26608^2&5^{15}17\cdot53093^{2}&2^{11}3^{4}101^{4}29221&-1&2&6.36599&21.9655\\
24304^2&3^{16}103^{3}127&73&7519&1&4.43014&12.3974\\
24064^2&2^{26}5^{12}1873&13^{10}37^{2}&2&1&14.8618&35.9005\\
22816^2&2^{30}5^{2}127\cdot353^{2}&1&6&1&8.19496&27.2210\\
22076^2&13^{2}251^{6}&2^{19}367^{3}&13&3&8.26382&58.6943\\
17264^2&3^{27}13^{4}&7^{3}29^{5}151^{2}&755&3&1.16852&14.8598\\
15888^2&5^{9}139^{6}&2^{2}11&5&1&15.1187&15.3391\\
15384^2&3^{22}13\cdot47^{2}263&2^{7}23^{8}&437&2&14.0272&18.8405\\
12664^2&19^{3}139^{6}&3\cdot13^{2}1049&13&3&5.02615&19.2050\\
12372^2&7\cdot131^{7}1373&3^{4}23^{6}1013^{2}&5&1&10.9177&16.7405\\
\hline\end{array}\end{gather*}

\newpage

{\bf Table 4}: Curves with $a=44$ and $G>4$.
\begin{gather*}\begin{array}{|c|c|c|c|c|c|c|}\hline|\Sha|&c&a&q&k&L&G\\\hline
479144^2&5^{9}139^{6}&44&104945&1&37.9640&96.2939\\
439312^2&5^{9}139^{6}&44&114998&1&30.4875&54.6806\\
321584^2&5^{9}139^{6}&44&5282&3&38.1136&31.3650\\
240848^2&5^{9}139^{6}&44&807&3&27.3470&23.9241\\
211784^2&5^{9}139^{6}&44&272857&1&9.19960&11.6672\\
183200^2&5^{9}139^{6}&44&92883&3&14.7483&7.99828\\
176832^2&5^{9}139^{6}&44&198770&1&15.0289&6.73874\\
144752^2&5^{9}139^{6}&44&243265&1&6.82733&5.77241\\
117552^2&5^{9}139^{6}&44&129965&1&4.10674&5.20828\\
117168^2&5^{9}139^{6}&44&5838&3&4.81256&6.52483\\
83536^2&5^{9}139^{6}&44&4865&1&2.67976&5.13812\\
83168^2&5^{9}139^{6}&44&1390&1&4.96930&17.8253\\
78048^2&5^{9}139^{6}&44&2310&3&6.78950&4.60258\\
70400^2&5^{9}139^{6}&44&302&1&7.63893&27.4015\\
55336^2&5^{9}139^{6}&44&3322&1&14.2300&5.10444\\
51872^2&5^{9}139^{6}&44&3315&3&2.50348&8.98022\\
49388^2&5^{9}139^{6}&44&4983&3&1.85105&6.63988\\
38920^2&5^{9}139^{6}&44&170&3&6.22360&11.1623\\
24428^2&5^{9}139^{6}&44&663&3&1.24148&4.45330\\
15888^2&5^{9}139^{6}&44&5&1&15.1187&15.3391\\
\hline\end{array}\end{gather*}

\section{Examples and strategy}

{\bf Example 1:} Consider the curve with $|\Sha|=479144^2$ in the first entry of Table 4.
The large value $G>96$ results from three circumstances. First, the factorization of
\begin{equation}
b=c-a=5^9\times139^6-44=3^2\times13^{10}\times17\times151\times4423
\end{equation}
is reasonably smooth, giving the radical of $abc$ as
\begin{equation}
r= 2\times3\times5\times11\times13\times17\times139\times151\times4423=6770408926710,\label{r44}
\end{equation}
which is significantly smaller than $c$, with $\lambda=\log(c)/\log(r)=1.49243$.
Secondly, $q=5\times139\times151=104945$ is a divisor of $r$. In general, the conductor
is given by
\begin{equation}
N=\frac{2^{s-1}q^2r}{{\rm gcd}(q,r)}\label{Nis}
\end{equation}
with an integer $s\in[0,5]$. In this example, $N=8qr$ and $q$ neatly disappears from
\begin{equation}\frac{G}{L}\approx\frac{r^{(\lambda-1)/2}}{2^{15/2}\pi}=2.53645,\label{G96}
\end{equation}
with an arithmetic-geometric mean in~(\ref{agm}) differing from unity  by less than $10^{-18}$
and the Tamagawa factors giving $C_1=64$. Finally, the large value $L=37.9640$ gives $G=96.2939$,
exceeding the previous record in~\cite{Nitaj}, with $|\Sha|=1832^2$ and $G=42.2653$.

{\bf Example 2:} Similar remarks apply to the second entry of Table 4, with $|\Sha|=439312^2$ and $G>54$.
Here $C_1=64$ and $q=2\times13\times4423=114998$ divides $r$ in~(\ref{r44}). The conductor in~(\ref{Nis}) is $N=16qr$, giving
$G/L\approx r^{(\lambda-1)/2}/(2^8\pi)=1.79354.$ Then $L=30.4875$ gives $G=54.6806$.

\newpage

The computation of central $L$-values proceeds by truncation of the infinite sum in~\cite{GZ}
\begin{equation}
L=L(1)=2\sum_{n>0}\frac{a(n)}{n} \exp(-2\pi n/\sqrt{N})\label{Lsum}
\end{equation}
where $a(n)=\prod_i a(p_i^{k_i})$ is determined by the factorization of $n=\prod_i p_i^{k_i}$ in powers of primes
and $a(1)=1$. For powers $k>1$ of primes that do not divide $N$, we have $a(p^k)=a(p^{k-1})a(p)-a(p^{k-2})p$.
If $p|N$, then $a(p^k)=a(p)^k$. If $p^2|N$, then $a(p)=0$.

Suppose that we truncate~(\ref{Lsum}) at $n=m$. Then $L$ is in error by a term
\begin{equation}
\delta_m=2\sum_{n>m}\frac{a(n)}{n} \exp(-2\pi n/\sqrt{N}).\label{dm}
\end{equation}
For $2\pi m\gg\sqrt{N}\gg1$, we expect that $\exp(2\pi m/\sqrt{N})\delta_m=O(1)$,
since $a(n)=O(n^{1/2})$  varies in sign for $n>m\gg1$. We can tolerate such an error
because $|\Sha|$ is assumed to be the square of an integer. Moreover, the maximal
value of $|\Sha|$ in an isogeny class is divisible by $4^t$ with an integer $t\in[0,5]$ determined
by the ratio of the largest and smallest values of the coefficients $C_k$ in~(\ref{agm}). 
If we assume that $\exp(2\pi m/\sqrt{N})|\delta_m|<K$, then we may truncate at
\begin{equation}
m\approx\frac{|\Sha|}{2\pi G}\log\left(\frac{K\sqrt{|\Sha|}}{2^tL}\right)\label{mK}
\end{equation}
where $|\Sha|=G\sqrt{N}$ is the maximal order of the Tate-Shafarevich group in the isogeny class.
In~\cite{GZ} it is stated (without proof) that we may set $K=4$, for sufficiently large $N$. This value
was adopted by Abderrahmane Nitaj in~\cite{Nitaj} and by Andrzej Dabrowski, in collaboration
with Mariusz Wodzicki in~\cite{DW} and with Lucjan Szymaszkiewicz in~\cite{DS}.

\newpage

The main {\em burden} in computing $|\Sha|$ comes from obtaining the coefficients $a(p)$ in~(\ref{Lsum})
for primes $p<m$, with $m$ estimated by setting $K=4$ in~(\ref{mK}). For this, I used the {\tt ellap} command of
Pari-GP~\cite{GP}, which implements the baby-step/giant-step method of Shanks and Mestre and runs in time roughly proportional
to $p^{1/4}$, implying that the total burden is roughly proportional to $m^{5/4}$. For each curve in the tables, 
I used a single thread of Pari-GP, running at 2.6~GHz on a modest laptop, with 3.2~GB of memory allocated
to this process. This allocation is sufficient  to allow storage of (and rapid access to)  values of $a(n)$ for $n<10^8$.
If $n$ is divisible by a prime power $p^k>10^8$, more work is required. For example,
with $m=3\times10^{10}$ and prime $p>10^8$ one might need to recompute $a(p)$ almost 300 times.
However, there is a saving that comes from the set ${\cal S}$ of primes $\ell$ such that $\ell^2|N$.
If $n$ is divisible by any prime $\ell\in{\cal S}$, then $a(n)=0$. Thus, for $|\Sha|>10^8$, I quantify the burden by
the integer nearest to
\begin{equation}
B=\left(\frac{|\Sha|}{(2\pi G)10^8}\log\left(\frac{4\sqrt{|\Sha|}}{2^tL}\right)\right)^{5/4}
\prod_{\ell^2|N}\left(1-\frac{1}{\ell}\right),\label{B}
\end{equation}
where~(\ref{mK}) has been divided by $10^8$, to make $B$ a convenient size.

For each of the 20 curves in Table 1, with $|\Sha|>250000^2$ and $G>12$, truncations of~(\ref{Lsum}) were monitored
as $n$ increased in steps of $10^7$ and the process was stopped when the approximate values of
$\sqrt{|\Sha|}/2^t$ for 20 successive steps differed from the same integer by less than 1/20. The wall-clock
time, in hours, was multiplied by 2.6 to give the actual cost $T$, in GHz-hours, rounded to the nearest integer.
I  found that~(\ref{B}) is a fair estimate of  the upper limit for $T$, with $T/B=0.8\pm0.2$.  Examples 1 and 2, with
$B=31$ and $B=61$,  took $T=28$ and $T=47$ GHz-hours, respectively.

\newpage

{\bf Example 3:} The first entry of Table~1 gives the curve with the record value
\begin{equation}
|\Sha|=2^6\times3^2\times13^2\times6211^2=1937832^2>3.755\times10^{12}\label{rec}
\end{equation}
achieved with  $G=153.084$. The factorization of
\begin{equation}
b=c-a=3^{19}\times11^{4}\times463^{5}-5^{4}\times19^{13}\times103=
2^{13}\times13^{9}\times29\times2441\times7673^{2}\label{brec}
\end{equation}
gives the radical of $abc$ as 
\begin{equation}
r=2\times3\times5\times11\times13\times19\times29\times103\times463\times2441\times7673=2111349279492568830
\end{equation}
for a triple with the good value  $\lambda=\log(c)/\log(r)=1.44935$, found by Frank Rubin in 2010~\cite{S}.
Its {\em merit}, as defined in~\cite{S}, is  $(\lambda-1)^2\log(r)\log(\log(r))=31.8833$, which is exceeded
by only 14 known triples.  A refinement of the $abc$ conjecture in~\cite{ABC} suggests an upper limit
of 48 for the merit of triples. The current record merit is less than 39. The quadratic twist $q=3\times5\times19=285$
divides $r$, giving a conductor $N=qr$. The Tamagawa factors give $C_3=576$ as the minimal integer in~(\ref{agm}),
where ${\rm AGM}(1,\sqrt{b/c})\approx0.999998$ differs little from unity. 
Thus we have $G/L\approx r^{(\lambda-1)/2}/(576\pi)>7.2366$
and hence $G>153$, thanks to a reasonably large central value $L=21.154$. It is this large value of $G$
that makes the burden~(\ref{B}) tolerable. With $B=958$ and $T/B\approx0.8$, I estimated
that it would take about 12 days at 2.6 GHz to determine  $|\Sha|$ running a single Pari-GP process, unthreaded. 
In the event, the elapsed time was 11.4 days, using the stopping criterion described above.

\newpage

{\bf Example 4:} The second entry of Table 1 gives a curve with $G>163$, which is the current record,
achieved
using a triple with $r=2249854771815161490$, $\lambda=1.41022$ and merit 26.6224,
found by Frank Rubin in 2019~\cite{S}. With $q=7$, $N=qr$ and $C_2=64$, this gives
$G/L=0.999898\,r^{(\lambda-1)/2}/(64\pi)=28.8963$ and hence $G=163.119$, with a rather modest
central value $L=5.64497$. The burden $B=189$ suggested a run time of about 2 days. In the event, it took 51 hours
to determine that $|\Sha|=804572^2$.

{\bf Example 5:} The third entry of Table 1 gives a curve with $G>162$, 
from a triple with $r=655240793892471930$, $\lambda=1.41235$ and merit 25.9073,
found by Tim Dokchitser in 2003~\cite{S}. With $q=23$, $N=qr$ and $C_1=160$, this gives
$G/L\approx r^{(\lambda-1)/2}/(160\pi)= 9.37574$ and hence $G=162.256$, with $L=17.3059$.
The burden is $B=167$. In the event, it took 39 hours to determine that $|\Sha|=793656^2$.

My strategy for selecting the 20 curves of Table 1 was based on experience
gained from the curves in Tables 2 and 3, where the burdens are smaller.
For a given triple $(a,b,c)$ and quadratic twist $q$, it takes little time to compute the 
arithmetic-geometric mean and Tamagawa factors that determine $G/L$.
The burden $B$ of determining an analytic value of $|\Sha|$ decreases
rather rapidly with increasing $G$, thanks to the factor of $1/G^{5/4}$ in~(\ref{B}).
The benefit of a large central value is much weaker, since  $L$ appears
only in the denominator of the logarithm of~(\ref{B}).  Hence I scanned
the 241 good triples in~\cite{S}, with square-free twists of magnitudes
$|q|<10^6$, selecting promising values of $G/L$ that might give
burdens $B<10^3$ for $L>1$. Only in such cases did I make a
very rough estimate of $L$, by truncating the sum in~(\ref{Lsum})
at $n=10^7$, which is far too small to give a reliable value of $L$.
If the very tentative estimate of $G$ seemed promising, I went up to $n=10^8$,
which takes about 6 minutes. 

\newpage

Thereafter it was usually clear whether the
curve was promising. However, there were some expensive disappointments,
where the estimate of $L$ decreased radically when going up to $n=10^9$,
or beyond, as well as happier cases, where it increased significantly.

In the 5 examples given so far, $q$ was a positive divisor of the radical $r$.
If $q$ is negative, the AGM in~(\ref{agm}) dilutes $G/L$, since $a<b$.
If $q$ is divisible by an odd prime $p$ coprime to $r$, then $G/L$ is diluted by a factor $p^{-1/2}$.
The following examples show that neither of these dilutions is an
insuperable obstacle to achieving $G>20$.
 
{\bf Example 6:}  The fourth curve in Table~1, with $|\Sha|=589080^2$, is based
on the triple of the record-breaking Example~3, but with a negative twist, $q=-95$, which dilutes
$G/L=0.215528\,r^{(\lambda-1)/2}/(256\pi)=3.50932$. However, it achieves
$G=24.5023$, thanks to $L=6.98206$.

{\bf Example 7:}  The fourth curve in Table~4, with $|\Sha|=240848^2$, is based
on the triple of Examples~1 and 2, but with a twist $q=3\times269=807$, 
where 269 is coprime to $r$ and hence dilutes
$G/L\approx r^{(\lambda-1)/2}/(32\sqrt{269}\pi)=0.874835$. However, it achieves
$G=23.9241$, thanks to $L=27.3470$.

\newpage

\section{Comparison with previous work}

In~\cite{W}, de Weger set the records $|\Sha|=224^2$ and $G=6.98260$, using an untwisted  triple  
from the relation
$7^3+5^{13}\times181=2^4\times3\times11\times13^2\times19^5$.

In~\cite{Nitaj}, Nitaj improved those records to $|\Sha|=1832^2$ and $G=42.2653$, twisting the
relation $5^{14}\times19+2^5\times3\times7^{13}=11^7\times37^2\times353$ by $q=11$.

Thereafter much work was done to find good $abc$ triples~\cite{S}, yet no improvement for $G$
was reported. Rather, Dabrowski and colleagues~\cite{DS,DW} determined large values of $|\Sha|$
using an Ansatz that produces rather low values of $G$ and hence entails burdens that are
much higher than in the present work. In~\cite{DW}, Dabrowski and Wodzicki
found 85 isogeny classes with $|\Sha|>10^6$, only two of which have $G>1$, while none has $G>2$.
The largest value $|\Sha|=63408^2$ has $G=0.911597$. 

More recently, Dabrowski and Szymaszkiewicz~\cite{DS} listed
40 isogeny classes with $|\Sha|>50000^2$, only three of which have  $G>1$, while none has $G>2$.
The following example illustrates how burdensome was the Ansatz that they adopted.
 
{\bf Example 8:} Let $a=29$, $b=3^{46}$ and 
$c=a+b$. The factorization $c=2\times17216879\times257390962660901$
shows that the radical of $abc$ is greater than $c$, with $r/c=3\times29=87$. With twist $q=3$,
we obtain a conductor $N=48r=4176c$ for the isogeny class (\ref{E1},\ref{E2},\ref{E3},\ref{E4}). 
The curves $E_1$ and $E_2$ give the minimal factors $C_1=C_2=2$ and hence the maximal value of $|\Sha|$.
The AGM in~(\ref{agm}) differs from unity by less than $10^{-21}$ and we obtain
$G/L\approx1/(8\sqrt{87}\pi)$, at very high precision. This is about 2200 times smaller than the value
of $G/L$ for the new record in Example~3.

\newpage 

In~\cite{DS} it was reported that Example 8 has $|\Sha|=1029212^2$, with $L=40.8169$ giving
the poor value $G=0.174117$. With $t=2$ in~(\ref{B}), from the maximal factor $C_3=32=4^2C_1$,
I obtain $B=578970$ as the notional burden. This is more than 600 times the burden $B=958$ for the new
record in Example~3.

It took less than 12 days, on a single thread running at 2.6 GHz with 3.2 GB of memory,
to increase the record for $|\Sha|$ by a factor of more than 3.5. It might have taken 20 years, with this 
modest resource, to check the previous record in~\cite{DS}, which was obtained by distributing
the task to a large number of processors.

\section{Statistics and prospects}

The strategy of Section~2 yielded 231 isogeny classes with $|\Sha|>10^8$ and $G>1$,
from twists of 78 good triples.
The scan of curves with $N<10^{18}$, for twists with $|q|<10^6$,
was reasonably thorough, yielding 173 classes. Yet it is inevitable that some cases were missed,
since the light initial screening by truncating~(\ref{Lsum}) at $n=10^7$ can be
misleading. For $N\in[10^{18},10^{19}]$, I was more selective, cherry-picking twists
that give good values of $G/L$, before spending the 6 minutes that were needed
to truncate at $n=10^8$, for a better estimate of $L$.  This resulted in 31 new classes. 
The remaining 27 classes, with $N>10^{19}$, were found by picking the very ripest cherries.
In particular, the 3 conspicuous cases with $|\Sha|>6\times10^{11}$ and $G>150$ 
emerged from candidates that were identified at an early stage of this work 
and then left to run their courses, while investigating less promising cases
on a second thread of the dual processor.

\newpage

In conclusion:
\begin{enumerate}
\item
Table~1 gives 20 classes with $\sqrt{|\Sha|}>250000$,  $G>12$ and burden $B<1000$.
\item
Table~2  gives 20 of the 29 classes with $\sqrt{|\Sha|}\in[100000,250000]$.
\item
Table~3  gives 20 of the 182 classes with $\sqrt{|\Sha|}\in[10000,100000]$.
\item
Table~4 gives 20 of the 38 classes with $c=5^{9}\times139^{6}$ and $a=44$.
\item
The triple with $c=2^{26}\times5^{12}\times1873$, $a=13^{10}\times37^2$ and merit 26.1686,
found by Tim Dokchitser in 2003~\cite{S}, is also rather fruitful. Its twists yielded 16 classes.
\item
The previous record prime $p$ dividing $|\Sha|$, found in~\cite{DS}, was $p=19861$. The list is hereby extended
to include the primes 25913, 27457, 36929 and 49253.
\item 
There are prospects of improving the record $|\Sha|=1937832^2\approx3.755\times10^{12}$, by
exploiting triples of high merit. By way of example, I was able to identify a class with
$|\Sha|=(4.99\pm0.01)\times10^{12}$  and  $G>50$.
\end{enumerate}

 {\bf Acknowledgements:} I thank Bart de Smit, for collating $abc$ triples of high merit, notably from 
Tim Dokchitser and Frank Rubin. Without the notable ABC triple of Bill Allombert, Karim Belabas and Henri Cohen, who
provide the excellent free-ware~\cite{GP} of Pari-GP, none of this work would have been feasible for me.
I am grateful to Kevin Acres and Michael St Clair Oakes for advice and encouragement.

 \newpage
 
\raggedright

}\end{document}